\documentclass[a4paper,10pt]{article}
\usepackage{amssymb,latexsym}
\usepackage{amsmath}
\usepackage{graphicx}

\newtheorem{theorem}{Theorem}
\newtheorem{lemma}{Lemma}

\def\qed{\ifhmode\unskip\nobreak\fi\quad\ifmmode\Box\else$\Box$\fi}

\title{Improvements on the density of maximal 1-planar graphs}
\author{J\'anos Bar\'at\\
\small  MTA-ELTE Geometric and Algebraic Combinatorics Research Group\\[-0.8ex]
\small \texttt{barat@cs.elte.hu}\\
and\\
G\'eza T\'oth\thanks{Supported by OTKA K-111827.}\\
\small Alfr\'ed R\'enyi Institute of Mathematics, 
\small Hungarian Academy of Sciences\\[-0.8ex]
\small \texttt{toth.geza@renyi.mta.hu}}

\begin{document}

\maketitle

\begin{abstract}
A graph is 1-planar if it can be drawn in the plane such that
each edge is crossed at most once. 
A graph, together with a 1-planar drawing is called 1-plane. 
Brandenburg et al. showed that 
there are maximal 1-planar graphs with only 
$\frac{45}{17}n + O(1)\approx 2.647n$ edges 
and maximal 1-plane graphs with only 
$\frac{7}{3}n+O(1)\approx 2.33n$ edges.
On the other hand, they showed that 
a maximal 1-planar graph has at least 
$\frac{28}{13}n-O(1)\approx 2.15n-O(1)$ edges, and
a maximal 1-plane graph has at least 
$2.1n-O(1)$ edges.

We improve both lower bounds to $\frac{20n}{9}\approx 2.22n$.
\end{abstract}

\section{Introduction}

In a {\em drawing} of a simple undirected graph $G$, 
vertices are represented by  distinct points in the plane
and edges are represented by simple
continuous curves connecting the corresponding points. 
For simplicity, the points and curves are also called vertices and edges, and
if it does not lead to confusion, 
we do denote them the same way
as the original vertices and edges of $G$, respectively. 
We assume that edges do not contain vertices in their interior. 
It follows from Euler's formula that a planar graph of $n$ vertices has
at most $3n-6$ edges. 
Also if a planar graph $G$ has less edges, then we can add edges to it so that the resulting graph has exactly
$3n-6$ edges and still planar. This holds even if we start with a fixed planar drawing of $G$. 
 
A drawing of a graph is {\it 1-planar}, if each edge is crossed
at most once. 
A graph is {\it 1-planar}, if it has a 1-planar drawing. 
It is {\it maximal} 1-planar, if we cannot add any edge to it so that the resulting graph is still 1-planar.
A graph together with a $1$-planar drawing is a {\it $1$-plane graph}.
It is {\it maximal} 1-plane, if we cannot add any edge to it so that the resulting drawing is still 1-plane.
The maximum number of edges of a 1-planar or 1-plane graph is $4n-8$ \cite{PT97}.
Recently, Brandenburg et al. \cite{ABG13,BEG13} observed a very interesting
phenomenon: there are maximal 1-planar and 1-plane graphs
with much fewer edges. 

\smallskip

\begin{theorem}[Brandenburg et al. $\cite{BEG13}$]
$(i)$ Let $e(n)$ be the minimum number of edges of a maximal $1$-planar graph 
with $n$ vertices. The following holds
$$\frac{28}{13}n-O(1)\approx 2.15n-O(1)\le e(n)\le \frac{45}{17}n + O(1)\approx 2.647n,$$

$(ii)$ Let $e'(n)$ be the minimum number of edges of a maximal $1$-plane graph 
with $n$ vertices. The following holds 
$$2.1n-O(1)\le e(n)\le \frac{7}{3}n + O(1)\approx 2.33n.$$
\end{theorem}

In this note, we improve both lower bounds. 

\begin{theorem} \label{main}
A maximal $1$-planar or $1$-plane graph has 
at least $\frac{20}{9}n-O(1)\approx 2.22n$ edges.
\end{theorem}

That is, $e(n), e'(n)\ge \frac{20}{9}n-O(1)$. 


\section{Preliminaries}

Our method is based on the ideas of Brandenburg et al. \cite{BEG13}. 
We also point out an error in \cite{BEG13}, but with our approach their proof goes through as well.
The following observations are essentially from their paper. 
We include the proofs for completeness.
Throughout this section, $G$ is a maximal $1$-plane graph. 
The edges of $G$ divide the plane into {\em faces}.
A face is bounded by edges and edge segments. 
These edges and edge segments end in vertices or crossings. 

\begin{lemma}\label{2v} 
{\rm (i)} There are at least two vertices on the boundary of each face.\\
{\rm (ii)} If $u$ and $v$ are two vertices on the boundary of a face, then they are adjacent.
\end{lemma}

\noindent {\bf Proof.} (i) Each face is bounded by at least three 
edges or edge segments, and has at least three vertices or crossings on its
boundary. 
Since there is at most one crossing on each edge, each edge segment contains 
a vertex as an endpoint. Therefore, there must be at least two vertices on the
boundary of the face.\\
(ii) 
Suppose that there are two vertices, $u$ and $v$ on the boundary of a face. 
Now $u$ and $v$ could be connected by a curve in the face without creating any
crossing.
Therefore, by the maximality of $G$, $u$ and $v$ are already connected. 
$\Box$ 

\begin{lemma} \label{nod1}
There are neither isolated vertices nor vertices of degree $1$ in $G$.
\end{lemma}

\noindent {\bf Proof.} Suppose that $v$ is an isolated vertex or a vertex of
degree 1 in face $F$ of $G$. 
Now $G\setminus\{ v\}$ is also maximal 1-planar, since if we can add
an edge to  $G\setminus\{ v\}$, then we could have added it to $G$. 
Therefore, by Lemma~\ref{2v}, $F$ has at least two vertices on its boundary,
different from $v$, and $v$ is adjacent to both of them, a contradiction. $\Box$

\begin{lemma} \label{ck4}
If $ab$ and $cd$ are crossing edges in $G$, then $a,b,c,d$ span a $K_4$ in $G$.
\end{lemma}

\noindent {\bf Proof.} 
Let $x$ be the crossing of $ab$ and $cd$. Since there are no other crossings
on $ab$ and $cd$, there is a face bounded by $ax$ and $xc$.
Now $a$ and $c$ are adjacent by Lemma~\ref{2v}. 
Similarly $a$ and $d$, $b$ and $c$, $b$ and $d$ are also adjacent. $\Box$ 

\smallskip

The smallest degree in $G$ is at least two by Lemma~\ref{nod1}. 
Following  \cite{BEG13}, we call vertices of degree two {\em hermits}.

\begin{lemma}
If a vertex $h$ has only two neighbors in $G$, say $u$ and $v$, then\\
{\rm (i)} $hu$ and $hv$ are not crossed by any edge,\\
{\rm (ii)}  $u$ and $v$ are adjacent in $G$.
\end{lemma}

\noindent {\bf Proof.} 
(i) Suppose to the contrary that $hu$ is crossed by an edge. 
By Lemma~\ref{ck4}, vertex $h$ has degree at least 3, a contradiction.\\
(ii) Since the only neighbors of $h$ are $u$ and $v$, and
edges $hu$ and $hv$ are not crossed, there is a face that has $u,x$ and
$v$ on its boundary. 
Therefore, $u$ and $v$ are adjacent by Lemma~\ref{2v}. $\Box$ 

\begin{lemma} 
Suppose that $h$ is a hermit, and its neighbors are $u$ and $v$. 
Delete $h$, $hu$, $hv$, $uv$ from $G$, and let $G'$ be
the resulting graph with the original embedding.
Let $F$ be the face of $G'$ that contains the point corresponding to vertex $h$.
Then $F$ has only two vertices on its boundary, $u$ and $v$.
\end{lemma}

\noindent {\bf Proof.} If there is another vertex on the boundary of $F$, then
we could connect it to $h$. Either without any crossing or with exactly one
crossing with edge $uv$ contradicting the maximality of $G$. $\Box$ 

\begin{figure}[h]
 \begin{center}
  \includegraphics[scale=0.4]{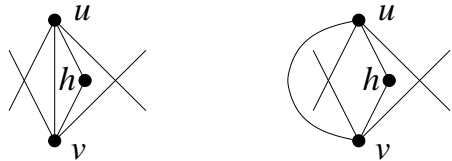}
 \end{center}
  \caption{Hermit $h$, surrounded by two pairs of crossing
edges.}
\label{hermit}
\end{figure}

We conclude that a hermit is surrounded by two pairs of crossing
edges, see Figure~\ref{hermit}. 

Remove all hermits from $G$. The resulting graph $\hat{G}$ with the
inherited drawing is the {\it skeleton} of $G$. 
Notice that $\hat{G}$ is also maximal 1-planar and each  vertex of $\hat{G}$ has degree at least $3$.


\subsection{A Correction}

In \cite{BEG13} the lower bound proofs rely on the following statement:

\smallskip

\noindent {\bf Claim.} \cite{BEG13} {\em 
Every edge of $\hat{G}$ is covered by a $K_4$ in $\hat{G}$.}

\smallskip

\begin{figure}[h]
 \begin{center}
  \includegraphics[scale=0.32]{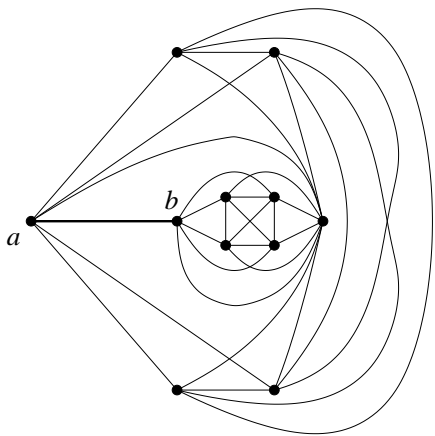}
 \end{center}
  \caption{A maximal 1-plane graph, $ab$ is not part of a $K_4$.}
\label{exceptional-kieg}
\end{figure}

However, this claim does not hold! 
See Figure~\ref{exceptional-kieg} for a counterexample.
Call an edge of  $\hat{G}$ {\em exceptional} if it is not part of a $K_4$ in $\hat{G}$. 
We have to deal with exceptional edges as well.

\begin{lemma} \label{exce}
Suppose that edge $ab$ of  $\hat{G}$ is exceptional. 
That is, $ab$ is not part of  a $K_4$ in $\hat{G}$. 
Let $F_1$ and $F_2$ be the faces bounded by $ab$.
Then

\noindent {\rm (i)} $F_1\neq F_2$, 

\noindent {\rm (ii)} $F_i$ has exactly three vertices on its
boundary $a,b$ and $f_i$ for $i=1,2$, and $f_1=f_2$,

\noindent {\rm (iii)} both $af_i$ and $bf_i$ are non-exceptional edges of  $\hat{G}$.
\end{lemma}

\begin{figure}[h]
 \begin{center}
  \includegraphics[scale=0.32]{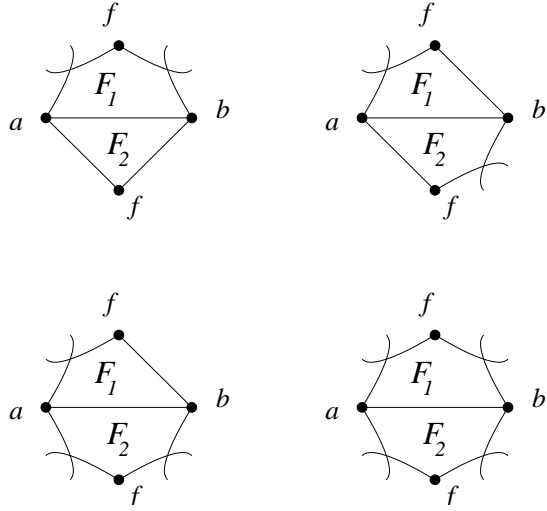}
 \end{center}
  \caption{The four possible types of exceptional edges $ab$}
\label{exceptional}
\end{figure}

\noindent {\bf Proof.} (i) Suppose that $ab$ is an exceptional edge of $\hat{G}$.
If edge $ab$ is crossed by another edge, then it is
part of a $K_4$ by Lemma~\ref{ck4}. 
Therefore, $ab$ does not participate in a crossing. 
Let $F_1$ and $F_2$ be the faces bounded by $ab$. 
If $F_1=F_2$, then $ab$ is a cut edge. 
In this case, by Lemma~\ref{2v}, both components have at least one other
vertex on the boundary of  $F_1=F_2$, and they can be connected. 
This contradicts the maximality of $\hat{G}$. 
Consequently, $F_1\neq F_2$. \\
(ii) If there is an edge from $a$ and an edge from $b$ which cross, then $ab$ is 
part of a $K_4$ by Lemma~\ref{ck4}. 
If there are at least four different vertices on the boundary of $F_i$, say
$a,b,x$ and $y$, then they form a $K_4$ by Lemma~\ref{2v}. 
We conclude that if $ab$ is exceptional, then $F_i$ has exactly three vertices on its
boundary, $a,b$ and $f_i$. 
If $f_1\neq f_2$, then we can connect them through $F_1$ and $F_2$, a contradiction again.
Therefore, $f_1=f_2$ and we denote it by $f$ for the rest of the proof.\\
(iii) Vertices $a,b$ and $f$ divide the boundary of $F_1$ into three parts. 
Between $a$ and $b$, we have edge $ab$ by assumption.
Between $a$ and $f$, we either have edge $af$, or two segments of edges. 
In the latter case, $af$ is an edge of  $\hat{G}$ and part of a $K_4$ by Lemma~\ref{ck4}.
We can argue similarly for face $F_2$. 
We conclude that $af$ is an edge of  $\hat{G}$ and part of a $K_4$, 
unless $af$ is on the boundary of both  $F_1$ and $F_2$. 
Now the degree of $a$ would be 2 in $\hat{G}$, which is impossible.
We can argue the same way for edge $bf$. $\Box$

\medskip

\noindent {\bf Remark.}
Each of the drawings on Figure~\ref{exceptional} can be extended to a
maximal 1-plane graph so that $ab$ is not part of a $K_4$. 
See Figure~\ref{exceptional-kieg} for an example. 


\section{Improvement of the lower bounds -- Proof of Theorem~\ref{main}}

Let $\hat{G}$ be the skeleton of a maximal 1-plane graph.
Recall that the skeleton $\hat{G}$ arises by removing each hermit from $G$ together with its two incident edges. 
The skeleton inherits its drawing from $G$, it is maximal 1-plane and each
vertex has degree at least three.
We distinguish three types of edges in $\hat{G}$: crossing, plain and exceptional. 
Edges that participate in a crossing are {\em crossing  edges}.
A crossing-free edge that is part of a $K_4$, is a {\em plain} edge.
Any other edge is {\em exceptional}.
Those edges are crossing-free and do not belong to a $K_4$.
Let $n({H})$, $c(H)$, $p(H)$, $e(H)$ denote the number of
vertices, crossing edges, plain edges, and exceptional edges of a graph $H$.
In particular, let $n=n(\hat{G})$, $c=c(\hat{G})$, $p=p(\hat{G})$, $e=e(\hat{G})$.
We prove the following crucial inequality involving these quantities.

\begin{lemma} \label{ineq}
 If $\hat{G}$ is the skeleton of any drawing  of a maximal $1$-planar graph $G$
 and $n\ge 4$, then 
\begin{equation}\label{eq1}
 9p+10e+7c\ge 20n-30.
\end{equation}
\end{lemma}

\noindent {\bf Proof.} 
We use induction on the pair $(e,n)$, ordered lexicographically.
If there is an exceptional edge, then we use the induction hypothesis on graphs with smaller $e$. 
If $e=0$, then we use induction on $n$.

Suppose that there is an exceptional edge $ab$ in $\hat{G}$.
Let $F_1$ and $F_2$ be the two faces bounded by $ab$.
By Lemma~\ref{exce}, $F_1\neq F_2$, and both $F_1$ and $F_2$ have exactly
three vertices on their boundaries, $a$, $b$, and $f$, see Figure~\ref{exceptional}.
The closure of $F_1\cup F_2$ divides the plane into two parts, say $S_1, S_2$.
Now  $\overline{S_i}$, the closure of $S_i$, intersects $\hat{G}$ in $G_i$ for $i=1,2$.
Remove the edge $ab$ and the interior of $F_1, F_2$ from $\hat{G}$.
Now two almost disjoint subgraphs $G_1$ and $G_2$ arise such that they have exactly one vertex $f$ in common,
$a\in G_1$ and $b\in G_2$. Both $G_1$ and $G_2$ are maximal 1-plane and both
have at least four vertices. 
Therefore, we can use the induction hypothesis on $G_1$ and $G_2$. 
For $i=1,2$, let  $n_i,c_i,p_i,e_i$ denote the number of
vertices, crossing edges, plain edges, and exceptional edges of $G_i$.
Now $9p_1+10e_1+7c_1\ge 20n_1-30$ and $9p_2+10e_2+7c_2\ge 20n_2-30$, 
where $e_1+e_2+1=e$, $n_1+n_2-1=n$ and $p_1+p_2=p$, $c_1+c_2=c$.
Therefore,
$9p+10e-10+7c\ge 20n+20-60$, and the statement follows.

We may now assume $e=0$, and we should prove $9p+7c\ge 20n-30$, where $n\ge 4$.
In what follows, we define an increasing sequence of subgraphs $G_0\subset G_1\subset\cdots\subset \hat{G}$
recursively and keep track of the number of vertices and edges of $G_i$. 
In every step, we maintain the inequality $9p+7c\ge 20n-30$.

Since there are no exceptional edges now, we can use the idea of Brandenburg et al. \cite{BEG13}.
They defined the {\it $K_4$-network} of $G$, which is an auxiliary graph ${\mathcal K}$.
Its vertex set corresponds to the $K_4$ subgraphs of ${\hat G}$.
Two vertices in ${\mathcal K}$ are adjacent if the corresponding subgraphs in ${\hat G}$ share a vertex.
Since $\hat{G}$ is connected and every edge is contained in a $K_4$, the graph ${\mathcal K}$ is connected.
Brandenburg et al. proved a lower bound on the number of edges of ${\hat G}$
by finding a certain spanning tree of ${\mathcal K}$ by an algorithm and 
investigating the number of edges of $\hat{G}$ involved in each step of the algorithm.

We go in their footsteps, but take a closer look.
We use a slightly more complex algorithm that sweeps through ${\hat G}$ rather than ${\mathcal K}$.


Let $G_0$ be a $K_4$ subgraph of $\hat{G}$.
Suppose that we have already defined $G_{i-1}$, a connected subgraph of $\hat{G}$, and now
we construct $G_i$. 
Therefore, the vertices, edges, subgraphs of $G_{i-1}$ are {\em old} and the ones of $G_i$ are {\em new}. 
We assume $9p(G_{i-1})+7c(G_{i-1})\ge
20n(G_{i-1})-30$.  This clearly holds for $i=1$.
To construct $G_i$ from $G_{i-1}$, we use one of the following operations in
this order of preference. 

\begin{enumerate}

\item Adding an edge between two old vertices.

\item  Adding a new vertex $x$ and all $K_4$'s spanned by $x$ and three old vertices.

\item Adding two new vertices, $x$ and $y$ and all $K_4$'s spanned by $x$, $y$
  and two old vertices.

\item Adding two new $K_4$'s such that they share a new vertex and each of them has a vertex in common with the old subgraph.

\end{enumerate}

If none of these operations can be executed, then let 
$G_{\mbox{\scriptsize final}}=G_{i-1}$ 
and the algorithm terminates.

Observe that $K_4$ has two different 1-planar drawings. 
Either all edges are crossing-free, or there is exactly one crossing.

We show that the output of the algorithm $G_{\mbox{\scriptsize final}}$ satisfies $G_{\mbox{\scriptsize final}}=\hat{G}$. 
We also show by induction, that $9p(G_i)+7c(G_i)\ge 20n(G_i)-30$ for every $i$. 
This is certainly true for $i=0$. 
Suppose $9p(G_{i-1})+7c(G_{i-1})\ge 20n(G_{i-1})-30$ for some $i$,
and now we construct $G_i$.
If we apply operation 1, then the left side of inequality (\ref{eq1})
increases by at least 5 (it is the case when we add an edge that crosses a
previously plain edge), while the right side does not change, so (\ref{eq1}) holds
for $G_i$ as well.

\begin{figure}[h]
 \begin{center}
  \includegraphics[scale=0.4]{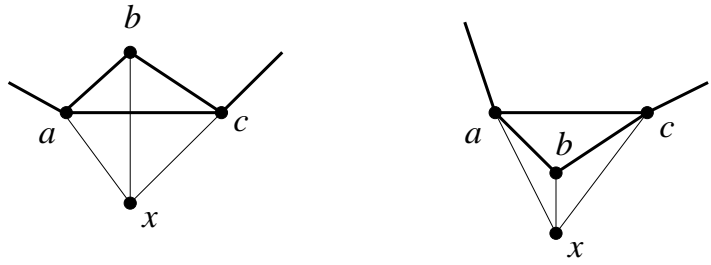}
 \end{center}
  \caption{Operation 2.}
\label{oper2}
\end{figure}

Suppose that we executed the second operation and we added exactly one $K_4$
with new vertex $x$ and old vertices 
$a$, $b$, $c$, see Figure~\ref{oper2}. 
Now either $xa$, $xb$ and $xc$ are plain edges in $G_i$, or
one of them, say $xc$,  crosses $ab$. The other two edges, $xa$ and $xb$ do not
cross an old edge since in this case there would be more $K_4$'s involving $x$. 
The left side of the inequality (\ref{eq1}) increased by $27$ or $23$ while the
right side increased by $20$. 
If we added more than one $K_4$, then we had to
add at least four edges adjacent to $x$. The addition of an edge increases the
left side of the inequality by at least $5$, so it increased by at least
$20$.

\begin{figure}[h]
 \begin{center}
  \includegraphics[scale=0.4]{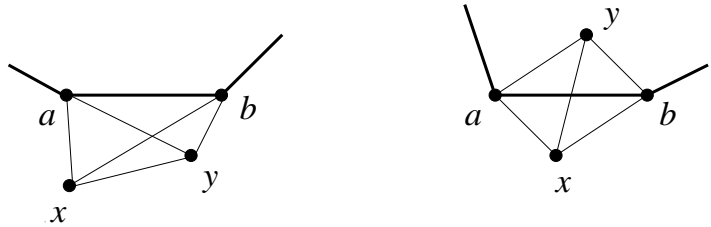}
 \end{center}
  \caption{Operation 3.} 
\label{oper3}
\end{figure}

Suppose that we executed the third operation, see Figure~\ref{oper3}.
Let $x$ and $y$ be the new vertices,
$a$ and $b$ the old vertices of a new $K_4$.
Edges $xa$, $xb$, $ya$, $yb$ cannot cross an old edge, since in that
case we find a $K_4$ with exactly three old vertices contradicting the preference order
of the operations. 
If $xy$ is not crossed by an old edge, then 
the left side of the inequality (\ref{eq1}) increases by $41$ or $45$, while the
right side increases by $40$. 
Suppose that the edge $xy$ crosses an old edge $cd$. 
Now $x, y, c, d$ form another $K_4$, and again none of the other
new edges crosses an old edge. 
We added at least eight new edges, so the
left-hand side increased by at least $40$ again.

\begin{figure}[h]
 \begin{center}
  \includegraphics[scale=0.4]{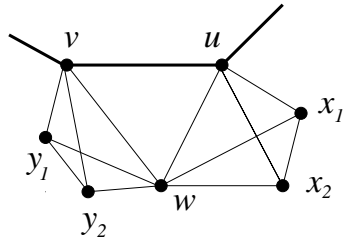}
 \end{center}
  \caption{Operation 4.}
\label{oper4}
\end{figure}

Suppose now that we arrive to a stage, 
where we cannot use any of the first three operations. 
Therefore, there is no $K_4$ in $\hat{G}$ that has exactly two or three old vertices.
Let $u$ be an old vertex that has at least one 
neighbor not in $G_{i-1}$.
Since $\hat{G}$ is connected, there is such a vertex. 
The graph $G_{i-1}$ is also connected, so $u$ has a neighbor in $G_{i-1}$ as well.
Order all neighbors of $u$ in the circular order the corresponding 
edges emanate from $u$. 

Recall that edges with a common endpoint do not cross.
Let $v$ and $w$ be consecutive neighbors of $u$ such that $v\in G_{i-1}$ and
$w\not\in G_{i-1}$. 
We distinguish four cases. 

Case 1: Both $uv$ and $uw$ are plain edges in $\hat{G}$.
We consider a $K_4$ that contains the edge $uw$. 
By the assumptions, this $K_4$ has exactly one vertex ($u$) in
$G_{i-1}$ and three vertices, say $w$, $x_1$, and $x_2$ not in $G_{i-1}$. 
The vertices $v$ and $w$ can be connected along $uv$ and $uw$, so by the
maximality of $\hat{G}$, they are adjacent in $\hat{G}$. 
We consider a $K_4$ that contains the edge $vw$.
By the assumptions, this $K_4$ has exactly one vertex ($v$) in
$G_{i-1}$ and three vertices not in $G_{i-1}$: $w$, $y_1$, and $y_2$. 
By the assumptions, none of the new edges crosses an old edge. 
If $x_1$, $x_2$, $y_1$, $y_2$ are all different, then we added five new vertices so the right-hand
side of (\ref{eq1}) increased by $100$, see Figure~\ref{oper4}.
The left-hand side increased by at least $100$ since adding two crossing $K_4$'s means 4 crossing and 8 plain edges
and $4\cdot 7+8\cdot 9=100$.
If $x_i=y_j$ for some $i, j$, then the situation is even better, the
calculation is very similar.

Case 2: The edge $uv$ is plain in  $\hat{G}$ and $uw$ is crossing.
In this case, $uw$ is crossed by an edge $ab$ of $\hat{G}$. 
Now $a, b, u, w$ form a $K_4$, so $a$ and $b$ are not in $G_{i-1}$. 
Vertex $v$ can be connected to $a$ or $b$
along $uv$, $uw$, and $ab$. Suppose that it is $a$, so by the maximality of
$\hat{G}$, $v$ and $a$ are adjacent. 
Consider a $K_4$ that contains edge $va$ and the one with vertices $a, b, u, w$. 
The calculation is very similar to the previous case.

Case 3: The edge $uw$ is plain in  $\hat{G}$ and $uv$ is crossing.
Let $ab$ be the edge that crosses $uv$. Now $a, b, u, v$ form a $K_4$, and 
$u$, $v$ are old vertices, so $a$ and $b$ are also old vertices. 
Vertex $w$ can be connected to $a$ or $b$, say, $a$, 
along $uw$, $uv$, and $ab$. So again by the maximality of $\hat{G}$,
$w$ and $a$ are adjacent. 
Consider a $K_4$ that contains $wa$ and a $K_4$ that
contains $uw$. The calculation is similar to the previous cases.

Case 4:  Both $uv$ and $uw$ are crossing edges in $\hat{G}$.
This case is the combination of the previous two cases. 
Edge $uw$ is crossed by $ab$, 
$uv$ is crossed by $cd$.
Now one of $a$ and $b$, say, $a$, and one of $c$ and $d$, say, $c$, can be
connected along $ab$, $uv$, $uw$, $cd$, so they are connected. 
Take the $K_4$ formed by $u, w, c, d$, and a $K_4$ the contains $ac$. 
The calculation is the same again.

In summary, we proved that we can always apply one of the four operations in
our algorithm, so the algorithm terminates when
$G_{\mbox{\scriptsize final}}=\hat{G}$.
On the other hand, we also proved $9p(G_{i})+7c(G_{i})\ge 20n(G_{i})-30$ for every $i$.
Therefore $9p(\hat{G})+7c(\hat{G})\ge 20n(\hat{G})-30$. This concludes the proof of
Lemma~\ref{ineq}. $\Box$

\subsection{Proof of Theorem~\ref{main}}

Recall that $e(n)$ ($e'(n)$) 
is the minimum number of edges of a maximal 1-planar (1-plane) graph 
with $n$ vertices.
Since every maximal 1-planar graph with any 1-planar drawing is a maximal
1-plane graph, $e(n)\ge e'(n)$. Therefore, Theorem~\ref{main} follows immediately
from the next result.

\begin{theorem}\label{seged}
Every maximal $N$-vertex $1$-plane graph has at least
$\frac{20}{9}N-\frac{10}{3}$ edges, where $N\ge 4$.
\end{theorem}

\noindent {\bf Proof.} 
Let $G$ be a maximal 1-plane graph, 
$N$ and $E$ denote the number of
vertices and edges, and $h$ denotes the number of hermits. 
Let $\hat{G}$ be the skeleton of $G$ and let 
$n=n(\hat{G})$, $c=c(\hat{G})$, $p=p(\hat{G})$, $e=e(\hat{G})$ denote the number of
vertices, crossing edges, plain edges, and exceptional edges of $\hat{G}$.

Every hermit is surrounded by two pairs of crossing edges.
A crossing pair of edges can participate in four such surroundings, on the
four sides of the crossing.
This gives us $c\ge h$. 
On the other hand, for each exceptional edge, each of the two neighboring 
cells has a pair of crossing edges on its boundary, these two crossings cannot 
participate in a surrounding of a hermit in that direction.
This shows $c\ge e$ and $c-e\ge h$.
Now $N=n+h$, $E=p+e+c+2h$. 

Let us minimize 
$$F(p,e,c,h,n)=E-\frac{20}{9}N$$ 
under the conditions 
\begin{equation} \label{eq2}
c-e\ge h,
\end{equation}
\begin{equation} \label{eq3}
9p+10e+7c\ge 20n-30,
\end{equation}
and
\begin{equation} \label{eq4}
p, e, c, h, n \ge 0
\end{equation}

$F(p,e,c,h,n)=E-\frac{20}{9}N=p+e+c+2h-\frac{20}{9}n-\frac{20}{9}h=
p+e+c-\frac{20}{9}n-\frac{2}{9}h$.

First we apply the following transformation:\\
$e'=e-\frac{9}{10}\varepsilon$, 
$p'=p+\varepsilon$, 
$h'=h+\frac{9}{10}\varepsilon$, 
$c'=c$, $n'=n$. 

Notice that if conditions (\ref{eq2}) and (\ref{eq3}) hold for 
$(p,e,c,h,n)$, then they also hold for 
$(p',e',c',h',n')$.
On the other hand, 
$F(p',e',c',h',n')=F(p,e,c,h,n)-\frac{1}{10}\varepsilon$.
Therefore, the five-tuple $(p,e,c,h,n)$ that minimizes $F(p,e,c,h,n)$
under conditions (\ref{eq2}) and (\ref{eq3}) has $e=0$.

For parameter $h$, the only condition is that $c\ge h$. If $c>h$ and we
increase $h$, then $F(p,0,c,h,n)$ decreases, and the conditions still
hold. Therefore, we may assume $c=h$. 
Now we have to minimize $F(p,0,c,c,n)=
p+\frac{7}{9}c-\frac{20}{9}n$ under the condition
$9p+7c\ge 20n-30$. We get immediately that the minimum of $F(p,0,c,c,n)$ under
the conditions is  $-\frac{10}{3}$.
Consequently $E-\frac{20}{9}N\ge -\frac{10}{3}$.

Therefore, $E\ge \frac{20}{9}N-\frac{10}{3}$ for any maximal 1-planar drawing with $N\ge 4$ vertices and $E$ edges.
 $\Box$

\medskip

\noindent {\bf Remark.} 
We believe that our bound is far from optimal. 
If our bound was close to optimal, then for some maximal 1-plane graph 
we would have to use operation 4. in almost every step
of the algorithm described in the proof of Lemma \ref{ineq}.
However, this seems impossible. 


\end{document}